\documentclass[12pt]{amsart}
\usepackage{amscd,amsmath,amsthm,amssymb,enumerate}
\usepackage[left]{lineno}
\usepackage{pstcol,pst-plot,pst-3d}
\usepackage{color}
\usepackage{pstricks}
\usepackage{stmaryrd}
\usepackage[utf8]{inputenc}
\usepackage{pstricks-add}
\usepackage[dvips]{graphicx}
\usepackage{epstopdf}
\usepackage{graphicx}
\newpsstyle{fatline}{linewidth=1.5pt}
\newpsstyle{fyp}{fillstyle=solid,fillcolor=verylight}
\definecolor{verylight}{gray}{0.97}
\definecolor{light}{gray}{0.9}
\definecolor{medium}{gray}{0.85}
\definecolor{dark}{gray}{0.6}

 \def\NZQ{\mathbb}               

 \def\ZZ{{\NZQ Z}}

 \def\frk{\mathfrak}               

 \def\mm{{\frk m}}

 \def\G{{\mathcal G}}

 \def\ab{{\mathbf a}}

 \def\cb{{\mathbf c}}

 \def\opn#1#2{\def#1{\operatorname{#2}}} 
 \opn\chara{char} \opn\length{\ell} \opn\pd{pd} \opn\rk{rk}
 \opn\projdim{proj\,dim} \opn\injdim{inj\,dim} \opn\rank{rank}
 \opn\depth{depth} \opn\grade{grade} \opn\height{height}
 \opn\embdim{emb\,dim} \opn\codim{codim}
 
 \opn\Tr{Tr} \opn\bigrank{big\,rank}
 \opn\superheight{superheight}\opn\lcm{lcm}
 \opn\trdeg{tr\,deg}
 \opn\reg{reg} \opn\lreg{lreg} \opn\ini{in} \opn\lpd{lpd}
 \opn\size{size} \opn\sdepth{sdepth}
 \opn\link{link}\opn\fdepth{fdepth}\opn\lex{lex}
 \opn\tr{tr}
  \opn\Hilb{Hilb}
 \opn\type{type}
 \opn\gap{gap}
 \opn\arithdeg{arith-deg}
 \opn\div{div} \opn\Div{Div} \opn\cl{cl} \opn\Cl{Cl}
 \opn\Spec{Spec} \opn\Supp{Supp} \opn\supp{supp} \opn\Sing{Sing}
 \opn\Ass{Ass} \opn\Min{Min}\opn\Mon{Mon}
 \opn\Ann{Ann} \opn\Rad{Rad} \opn\Soc{Soc}
 \opn\Im{Im} \opn\Ker{Ker} \opn\Coker{Coker} \opn\Am{Am}
 \opn\Hom{Hom} \opn\Tor{Tor} \opn\Ext{Ext} \opn\End{End}
 \opn\Aut{Aut} \opn\id{id}
 
 \opn\nat{nat}
 \opn\pff{pf}
 \opn\Pf{Pf} \opn\GL{GL} \opn\SL{SL} \opn\mod{mod} \opn\ord{ord}
 \opn\Gin{Gin} \opn\Hilb{Hilb}\opn\sort{sort}
 \opn\PF{PF}\opn\Ap{Ap}
 \opn\mult{mult}
 \opn\bight{bight}
 \opn\aff{aff}
 \opn\relint{relint} \opn\st{st}
 \opn\lk{lk} \opn\cn{cn} \opn\core{core} \opn\vol{vol}  \opn\inp{inp} \opn\nilpot{nilpot}
 \opn\link{link} \opn\star{star}\opn\lex{lex}\opn\set{set}
 \opn\width{wd}
 \opn\Fr{F}
 \opn\QF{QF}
 \opn\G{G}
 \opn\type{type}\opn\res{res}
 \opn\conv{conv}
 \opn\Shad{Shad}
 \opn\ln{ln}
 \opn\In{in}
 \opn\gr{gr}
 
 \def\pot#1#2{#1[\kern-0.28ex[#2]\kern-0.28ex]}

 \opn\dirlim{\underrightarrow{\lim}}
 \opn\inivlim{\underleftarrow{\lim}}
 \let\iso=\cong

 \def\Implies{\ifmmode\Longrightarrow \else
         \unskip${}\Longrightarrow{}$\ignorespaces\fi}
 \def\implies{\ifmmode\Rightarrow \else
         \unskip${}\Rightarrow{}$\ignorespaces\fi}
 \def\iff{\ifmmode\Longleftrightarrow \else
         \unskip${}\Longleftrightarrow{}$\ignorespaces\fi}

 \let\:=\colon
 \newtheorem{Theorem}{Theorem}[section]
 
 \newtheorem{Corollary}[Theorem]{Corollary}
 \newtheorem{Proposition}[Theorem]{Proposition}
 \newtheorem{Remark}[Theorem]{Remark}
 
 \newtheorem{Example}[Theorem]{Example}
 
 \newtheorem{Definition}[Theorem]{Definition}
 \newtheorem{Notation}[Theorem]{Notation}

 \let\epsilon\varepsilon
 \let\kappa=\varkappa
 \def\qed{\ifhmode\textqed\fi
       \ifmmode\ifinner\quad\qedsymbol\else\dispqed\fi\fi}
 \def\textqed{\unskip\nobreak\penalty50
        \hskip2em\hbox{}\nobreak\hfil\qedsymbol
        \parfillskip=0pt \finalhyphendemerits=0}
 \def\dispqed{\rlap{\qquad\qedsymbol}}
 
 \opn\dis{dis}
 \def\pnt{{\raise0.5mm\hbox{\large\bf.}}}
 
 \opn\Lex{Lex}

\begin{document}

\title {Toric ideals  which are determinantal}

\author {Reza Abdolmaleki and Rashid Zaare-Nahandi }

\address{Reza Abdolmaleki, Department of Mathematics, Institute for Advanced Studies in Basic Sciences (IASBS), 45195-1159, Zanjan, Iran}
\email{abdolmaleki@iasbs.ac.ir}

\address{Rashid Zaare-Nahandi, Department of Mathematics, Institute for Advanced Studies in Basic Sciences (IASBS), 45195-1159, Zanjan, Iran} 
\email{rashidzn@iasbs.ac.ir}

\dedicatory{ }

\begin{abstract}  

Given any equigenerated monomial ideal $I$ with the property that the defining ideal $J$ of the fiber cone $ F(I)$ of $I$ is generated by quadratic binomials, we introduce a matrix such that the set of its binomial $2$-minors is a generating set of  $J$. In this way, we characterize the fiber cone of sortable and Freiman ideals.
\end{abstract}

\thanks{}



\maketitle

\setcounter{tocdepth}{1}

\section*{Introduction}
Let $K$ be a field and $S = K[x_1, \ldots, x_n]$ the polynomial ring in variables $x_1,\ldots,x_n$  over $K$. For a graded ideal $I \subset S$  the fiber cone $F(I)$ of $I$ is the the standard graded $K$-algebra $\bigoplus_{k \geq 0} I^k/\mm I^k$, where $\mm$ denotes the unique maximal graded ideal of $S$. Indeed, $F(I)=\mathcal{R}(I)/\mm \mathcal{R}(I)$, where $\mathcal{R}(I) =\bigoplus_{k \geq 0} I^kt^k \subset S[t]$ is the Rees ring of $I$. 

 Let $I$ be a monomial  ideal with $ G(I)= \{u_1,  \ldots , u_q\}$ and $T=K[t_{u_1}, \ldots , t_{u_q}]$ be the polynomial ring in variables $t_{u_1}, \ldots , t_{u_q}$ over $K$.  The $K$-algebra homomorphism $T \rightarrow  F(I)$,  $t_{u_i} \mapsto u_i +\mm I$ induces the isomorphism $F(I) \iso T/J$. The ideal $J$ is called the defining ideal of $F(I)$. Finding the minimal generators of $J$  and the algebraic properties of $F(I)$ is a difficult problem even in concrete cases. For example, for the symmetric ideals 
\begin{eqnarray*}
&&I_1=(x_{1}^{11}, x_{1}^9x_{2}^2, x_{1}^7x_{2}^4,  x_{1}^{6}x_{2}^{5}, x_{1}^5x_{2}^6, x_{1}^4x_{2}^7, x_{1}^{2}x_{2}^{9},x_{2}^{11 }) ,\\
&&I_2=(x_{1}^{11}, x_{1}^{10}x_{2}, x_{1}^7x_{2}^4, x_{1}^{6}x_{2}^{5}, x_{1}^5x_{2}^6, x_{1}^4x_{2}^7, x_{1}x_{2}^{10}, x_{2}^{11 })
\end{eqnarray*}
 generated in degree $11$ in $K[x_1, x_2]$, one can check by CoCoA \cite{Co} that the defining ideal of $F(I_1)$ is an  ideal generated by quadratic binomials and  $F(I_1)$ is a Cohen-Macaulay algebra, while the minimal generating set of the defining ideal of $F(I_2)$ includes binomials in degrees $2$ and $4$,  and  $F(I_2)$ is not Cohen-Macaulay. Note that, $G(I_1)$  and $G(I_2)$  just differ in 2 monomials. We recall that a monomial ideal  $I \subset S = K[x_1, x_2]$ with $G(I) = {u_1, \ldots , u_q}$ and $u_i =x_1^{a_i}x_2^{b_i}$ satisfying the properties $a_1 > a_2 > \ldots > a_q = 0$ and $0 = b_1 < b_2 < \ldots < b_q$ is called
a {\em symmetric ideal}, if $b_i = a_{q-i+1}$ for $i = 1, \ldots, q$. The fiber cones of symmetric ideals with $4$ generators are well studied in \cite{HQS} and \cite{HZ}. Moreover,  \cite{HQS} includes a characterization of the fiber cones  of concave and convex monomial ideals  and their algebraic properties. 

In this paper,  for an equigenerated monomial ideal $I$ with the property that the defining ideal $J$ of  $F(I)$ is generated by quadratic binomials, we  interpret $J$ as the set of  binomial $2$-minors of a special matrix. By an equigenerated monomial ideal we mean an ideal generated by monomials in a single degree. In this paper we consider that a binomial has exactly two non-zero monomials.

In Section~\ref{1}, we associate to a sortable ideal $I$ a matrix $T_I$, such that the defining ideal $J$ of $F(I)$ is generated by the set of  binomial $2$-minors of $T_I$ (Theorem~\ref{main}). Moreover, we show that equigenerated $\cb$-bounded strongly stable monomial ideals, in particular Veronese type ideals, are sortable and hence, their fiber cones are Cohen-Macaulay normal domains and reduced Koszul algebras. 

 In Section~\ref{2},   for any equigenerated monomial ideal $I \in K[x_1, \ldots x_n]$ with $n \geq 3$,  we show that if the defining ideal $J$ of the fiber cone $F(I)$ is generated by quadratic binomials, then $J$ is generated by the set of  binomial $2$-minors of $T_I$ (Theorem~\ref{more}).  For the case $I \subset K[x_1, x_2]$,   we associate to $I$ a matrix $\mathcal{T_I}$ in a different way, and show that the defining ideal of $F(I)$ is generated by the set of  binomial $2$-minors of $\mathcal{T_I}$ (Theorem~\ref{tor}). In particular, we determine the fiber cone of Freiman ideals. A Freiman ideal  $I$  is an equigenerated monomial ideal such that $\mu (I^2)=\ell (I) \mu (I)- {\ell (I) \choose 2}$, where $\mu (I)$ is the minimal number of generators of $I$, and $\ell (I) $  denotes the analytic spread of $I$ which is by definition the Krull dimension of $F(I)$.  Freiman ideals are studied in \cite{HHZ} and \cite{HZ1}.

\section{The fiber cone of sortable ideals}
\label{1}

Let $S=K[x_1,\ldots,x_n]$ be the polynomial ring in the variables $x_1,\ldots,x_n$ over a field  $K$ . We denote by $\mm$ the unique maximal graded ideal of $S$.  Let $I$ be a monomial ideal of $S$ and $ G(I)= \{u_1,  \ldots , u_q\}$ be the minimal set of monomial generators of $I$. The fibre cone $F(I)$ of $I$ is defined as the standard graded $K$-algebra $\bigoplus_{k \geq 0} I^k/\mm I^k$. Let $T=K[t_{u_1}, \ldots , t_{u_q}]$, where $t_{u_1}, \ldots , t_{u_q}$ are independent variables. Then $F(I) \iso T/J$, where $J$ is the kernel of the $K$-algebra homomorphism $T \rightarrow  F(I)$ with $t_{u_i} \mapsto u_i +\mm I$. In this section we determine the fiber cone of  sortable ideals and introduce several classes of ideals which are sortable.

\medskip
Let $d$ be a positive integer and $S_d$ be the $K$-vector space generated by monomials of degree $d$ in $S$.  For monomials $ u, v \in S_d$ we write $uv = x_{i_1}x_{i_2} \ldots x_{i_{2d}}$ with $ 1 \leq  i_1 \leq i_2 \leq  \ldots \leq i_{2d} \leq n$.
The pair $(u',v')$ with $ u' = x_{i_1}x_{i_3} \ldots x_{i_{2d-1}}$ and $ v' = x_{i_2}x_{i_4} \ldots x_{i_{2d}}$ is called the {\em sorting} of $(u,v)$.  So we get the map
\[
\sort : S_d \times S_d \rightarrow S_d \times S_d, (u, v) \mapsto (u', v'),
\]
called {\em sorting operator}. A pair $(u,v)$ is called {\em sorted} if $\sort (u,v)=(u,v)$, otherwise it is called {\em unsorted}.
It is shown in \cite[Section~6.2]{EH} that the pair $(u,v)$ with $u=x_{i_1} x_{i_2} \ldots x_{i_d}$ and $v=x_{j_1} x_{j_2} \ldots x_{j_d}$  is sorted if and only if
\begin{eqnarray}
\label{sort}
i_1 \leq j_1 \leq i_2 \leq j_2 \leq \ldots \leq  i_d \leq j_d.
\end{eqnarray}
 Note that  if $(u,v)$ is sorted, then $u  \geq_{ \lex} v$,  where $\geq_{ \lex} $ denotes the lexicographic order on $\Mon(S)$, the set of monomials of $S$.
\begin{Definition}
{\em (a)  A set of monomials $A \subset S_d$  is called {\em sortable} if $ \sort(A\times A) \subset A \times A$.}

{\em (b) An equigenerated monomial ideal $I$ is called a {\em sortable ideal}, if  $G(I)$ is a sortable set.}
\end{Definition}
Let $A \subset S_d$ be a sortable set of monomials and $I$ be the ideal generated by $A$. We denote by $ K[A]$ the semigroup ring generated over $ K$ by $A$. Let $T = K[{t_u : u \in A}]$ be the polynomial ring with the order on variables given by $t_u > t_v$ if $u >_{lex} v$ . Also, let  $\varphi : T \rightarrow K[A]$  be the $K$-algebra homomorphism defined by $ t_u \mapsto u$ for all $ u \in A$ and $P_A$ be
the kernel of $\varphi$. Since the ideal $I$ is equigenerated, $F(I) \iso K[A]$ (see \cite[the proof of Corollary~1.2]{HMZ}) and so the toric ideal $P_A$ is the defining ideal $J$ of $F(I) $  in the representation $F(I) =T/J$ of the fiber cone of $I$.

The following well known theorem plays an important role in the proof of the main theorem of this section (See \cite[Theorem~6.16]{EH}).
\begin{Theorem}
\label{gr}
Let $K[A]$ be a $K$-algebra generated by a sortable set of monomials $A \subset S_d$ and $P_A\subset R$ its toric ideal. Then
\[
\mathcal{G} = \{ t_ut_v - t_u't_v' : u, v \in A, (u, v) \text{ unsorted} , (u', v') = \sort(u, v)\}
\]
is the reduced Gr\"obner basis of $P_A$ with respect to the sorting order.
\end{Theorem}
It follows from \cite[Theorem~6.15]{EH} that  the ideal $\In_<(\mathcal{G})$ is generated by the monomials $t_ut_v$ where $(u, v)$ is unsorted.

Before stating the main theorem of this section, we recall that an affine semigroup $H$  generated by the set $\mathcal{H} = \{h_1, \ldots, h_q\} \subset \ZZ^n$ is called {\em normal} if it satisfies the following condition: if $mg \in H$ for some $g \in \ZZ H$ and $m > 0$, then $g \in H$, where $\ZZ H$ is the subgroup of  $\ZZ^n$ generated by $\mathcal{H}$.  Also, a domain $R$ is called {\em normal} if it is integrally closed, that is $R=\bar{R}$ where $\bar{R}$ is the integral closure of $R$.
\begin{Notation}
\label{associate}
Let $I \subset S$ be an ideal  generated minimally  by a  set of monomials of degree $d$. One can consider $G(I)=\{u_1, \ldots, u_{q}\}$ as a subset of $G(\mm^d)$. Let $M$ be the matrix which  the entries of  its $i$-th row are the monomials of $G(\mm^d)$ containing $x_i$,  ordered lexicographically from left to right, for $i=1, \ldots, n$. This matrix has $n$ rows and ${ n+d-2 \choose d-1}$ columns. We replace the entries of $M$  belonging to $G(\mm^d) \setminus G(I)$  by $0$, remove its zero columns  and denote the obtained matrix by  $M_I$.  Then we replace any non-zero element $u$  of $M_I$  by indeterminate $t_u$ of $R=K[t_{u_1}, \ldots, t_{u_q}]$. Denote this matrix by $T_I$ and call it  the {\em matrix associated to $I$}.
\end{Notation}

\begin{Example}
\label{M}
Let $n=3$ and $d=3$. The matrix $M$ is the following:
\begin{align*}
M=\begin{pmatrix}
x_1^3      & x_1^2x_2  & x_1^2x_3  & x_1x_2^2 & x_1x_2x_3 & x_1x_3^2 \\ 
x_1^2x_2 & x_1x_2^2  & x_1x_2x_3 & x_2^3      & x_2^2x_3  & x_2x_3^2 \\ 
x_1^2x_3 & x_1x_2x_3 & x_1x_3^2  & x_2^2x_3 & x_2x_3^2  & x_3^3 
\end{pmatrix}.
\end{align*}
Let $I$ be the Veronese type ideal $I=(x_1^3, x_1^2x_2, x_1^2x_3, x_1x_2^2, x_1x_2x_3, x_2^2x_3) \subset K[x_1,x_2,x_3]$. We set $u_1=x_1^3$, $u_2=x_1^2x_2$, $u_3=x_1^2x_3$, $u_4=x_1x_2^2$, $u_5=x_1x_2x_3$ and $u_6=x_2^2x_3$. Let $F(I)=K[t_{u_1}, \ldots, t_{u_6}]/J$. We have
\begin{align*}
T_I=\begin{pmatrix}
 t_{u_1}& t_{u_2}  &  t_{u_3} &  t_{u_4}  &  t_{u_5} \\ 
 t_{u_2} &  t_{u_4}  & t_{u_5}  & 0 &  t_{u_6 }\\ 
 t_{u_3} &   t_{u_5} &  0& t_{u_6}  & 0 
\end{pmatrix}.
\end{align*}
\end{Example}

\begin{Theorem}
\label{main}
Let $I \subset S$ be a sortable ideal with $G(I)= \{u_1, \ldots , u_q \}$ and  the fiber cone $F(I)=K[t_{u_1}, \ldots, t_{u_q}]/J$. Also,  let $T_I$  be the matrix associated to $I$, introduced in Notation~\ref{associate}.
\item[(a)] The toric ideal $J$ is generated by the set of binomial $2$-minors of $T_I$. Indeed,
\[
J=(t_ut_v - t_{u'}t_{v' }: u, v, u' ,v' \in G(I),   \begin{pmatrix}
t_{u}& t_{u'}\\ 
t_{v'}& t_{v}\\ 
\end{pmatrix} \text{is a submatrix of} \ T_I).
\]
\item[(b)]  $F(I)$ is a reduced Koszul algebra.
\item[(c)] $F(I)$ is a Cohen-Macaulay normal domain.
\end{Theorem}
\begin{proof} 
(a) For $i=1, \ldots , n$, dividing the entries of the $i$-th row of the matrix $M$ by $x_i$, we get a matrix whose entries of all rows  are the monomials of $G(\mm^{d-1})$ ordered lexicographically from left to right. This implies that the set of binomial $2$-minors of $T_I$ includes in $J$.

Now, we show that for the  monomials $u, v, u', v' \in G(I)$ if the binomial $f=t_ut_v-t_{u'}t_{v'} $ belongs to $J$, then $f$ is a $2$-minor of $T_I$. By theorem~\ref{gr}
\[
\mathcal{G} = \{ t_ut_v - t_{u'}t_{v' }: u, v \in G(I),  (u, v) \text{unsorted}, (u', v') = \sort(u, v)\}
\]
is the reduced Gr\"obner basis of $J$ with respect to the sorting order. We show that if monomials $u,v \in G(I)$ are unsorted, then $u$  and $v$ are the $ ij$-th and the $ kl$-th  entries of the matrix $M_I$ respectively, such that $i \neq k$ and $j \neq l$  and  that $\sort(u,v)=(u',v')$, where $u'$  and $v'$ are the $ il$-th and the $ kj$-th  entries of the matrix $M_I$ respectively. This imples that   the determinats of $2\times2 $ submatrices of $T_I$ which have no zero entries, form a Gr\"obner basis of $J$. Let  $(u,v)$  with $u=x_{i_1} x_{i_2} \ldots x_{i_d}$ and $v=x_{j_1} x_{j_2} \ldots x_{j_d}$ be an unsorted pair in $G(I) \times G(I)$.  Notice that $u$ and $v$ belong to different columns of  $M$. Indeed, if  two different monomials $w,w'$ belong to the same column of $M$, we have $w=x_pw_1$ and $w'=x_qw_1$ for $1 \leq p\neq q \leq n$ and a  monomial $w_1 \in  G(\mm^{d-1})$ . So, by \eqref{sort} the pair $(w,w')$ is sorted. On the other hand, since $(u,v) $ is unsorted, $u$ is divisible by $x_i$ and $v$ is divisible by $x_j$  for some $1 \leq i\neq j \leq n$,  because otherwise $u=v=x_r^d$ for a variable $x_r$, a contradiction.  Hence,  we can find $u$ and $v$  in different rows of  the matrix $M$ (although they may appear  in the same row as well). So we assume that $u=u_{ij}$ and $v=v_{kl}$ such that $i \neq k$ and $j \neq l$. Since $G(I)$ is a sortable set, it follows that $\sort(u,v) \in G(I) \times G(I)$.  Assume that $\sort(u,v)=(u',v')$. So $t_ut_v - t_{u'}t_{v'} $ belongs to the reduced Gr\"obner basis of $J$ by Theorem~\ref{gr}. Note that $u_{ij} =x_iu_1$ and $v_{kj} = x_ku_1$ for a monomial $u_1 \in G(\mm^{d-1})$. Similarly, $u_{il} =x_iu_2$ and $v_{kl} = x_ku_2$ for a monomial $u_2 \in G(\mm^{d-1})$. Therefore, $u_{ij}v_{kl}=u_{il}v_{kj}=x_ix_ku_1u_2$ and hence  $\sort(u_{ij},v_{kl})=\sort(u_{il},v_{kj})$.  Suppose that $(u_{il},v_{kj})$ is unsorted. It follows from Theorem~\ref{gr} that $t_{u_{il}}t_{v_{kj}} - t_{u'}t_{v'}$ belongs to the reduced Gr\"obner basis of $J$ which is a contradiction, because $(t_ut_v - t_{u'}t_{v'})- (t_{u_{il}}t_{v_{kj}} - t_{u'}t_{v'})=t_ut_v -t_{u_{il}}t_{v_{kj}} $ belongs to the reduced Gr\"obner basis of $J$ (note that $t_ut_v - t_{u'}t_{v'} $ belongs to the reduced Gr\"obner basis of $J$). Therefore, $\sort(u_{ij},v_{kl})=(u_{il},v_{kj})$. So, the assertion follows from Theorem~\ref{gr}. Notice that this Gr\"obner basis  is not necessary reduced.

(b) It follows from \cite[Theorem~6.15]{EH}) that $\In_<(J)$  is a square-free monomial ideal. This yields that $J$ is a radical ideal (see \cite[Theorem~3.3.7]{HH})   and hence $F(I)$ is a reduced algebra. Moreover, since $J$ has a quadratic Gr\"obner basis,  $F(I)$ is Koszul  by a well known result of Fr\"oberg (see \cite[Theorem~6.7]{EH}).

(c) Since $\In_<(J)$ is a squarefree monomial ideal,  it follows from a result by Sturmfels (\cite[Proposition~13.15]{S})   that $F(I)$ is normal. Moreover, by a result of Hochster (\cite[Theorem~1]{Ho}) $F(I)$  is Cohen-Macaulay.
\end{proof}

In the rest of this section we show that  any equigenerated $\cb$-bounded  strongly stable monomial ideal is sortable.

 Let $\cb=(c_1,\ldots,c_n)$ be an integer vector with $c_i\geq 0$. The monomial  $u=x_1^{a_1}\cdots x_n^{a_n}$  is called {\em $\cb$-bounded}, if
 $\ab\leq \cb$, that is, $a_i\leq c_i$  for all $i$.
Let $I=(u_1,\ldots, u_m)$ be a monomial ideal.  We set
\[
I^{\leq \cb}=(u_i\:\; \text{ $u_i$  is  $\cb$-bounded}).
\]
We also set $m(u)=\max\{i\:\; a_i\neq 0\}$. The following definition is obtained from \cite{AHZ}.
\begin{Definition}
{\em Let $I\subset S$ be a  $\cb$-bounded monomial ideal.}
\begin{enumerate}
{\em \item[(a)]  $I$ is called {\em $\cb$-bounded  strongly stable} if for all  $u\in G(I)$ and all $i<j$  with $x_j|u$  and    $x_iu/x_{j}$ is $\cb$-bounded, it follows that $x_iu/x_{j}\in I$.}

{\em \item[(b)] $I$ is called {\em $\cb$-bounded  stable} if for all $u\in G(I)$ and all $i<m(u)$ for which   $x_iu/x_{m(u)}$ is $\cb$-bounded, it follows that $x_iu/x_{m(u)}\in I$.}
\end{enumerate}
\end{Definition} 
It is clear that a $\cb$-bounded  strongly stable monomial ideal is $\cb$-bounded  stable.

 The smallest $\cb$-bounded strongly stable ideal containing  $\cb$-bounded monomials $u_1,\ldots,u_m$ is denoted by $B^\cb(u_1,\ldots,u_m)$. A monomial ideal $I$ is called a {\em $\cb$-bounded  strongly  stable principal idea}l, if there exists a $\cb$-bounded monomial $u$ such that $I=B^{\cb}(u)$.  The smallest strongly stable ideal containing $u_1,\ldots,u_m$ (with no restrictions on the exponents) is denoted $B(u_1,\ldots,u_m)$. The monomials $u_1,\ldots,u_m$ are called {\em Borel generators}  of $I=B(u_1,\ldots,u_m)$.

\begin{Proposition}
Let $I=B^\cb(u_1,\ldots,u_m)$ be an equigenerated $\cb$-bounded strongly stable monomial ideal. Then $I$ is a sortable ideal.
\end{Proposition}
\begin{proof}
First we prove the assertion for the case $\cb$-bounded strongly stable principal ideal $I=B^\cb(u_k)$ where $1 \leq k \leq m$. Assume that $v,w \in G(B^\cb(u_k))$ and $\sort (v,w)=(v',w')$. For this purpose we first  show that $v', w'$ are $\cb$-bounded monomials. So, we must check that for all $i \in \{1, \ldots, n\}$, the degrees of $x_i$ in $v'$ and $w'$ are not greater than $c_i$. Let $\deg_{x_i}(v) = a_i$ and $\deg_{x_i}(w) = b_i$. Note that $a_i, b_i \leq c_i$. If $a_i+b_i$ is even, $\deg_{x_i}(v')= \deg_{x_i}(w') = (a_i+b_i)/2$ by the definition of the sorting operator. Therefore, $\deg_{x_i}(v'), \deg_{x_i}(w')  \leq c_i$. Now let $a_i+b_i$  be an  odd integer. Then $\deg_{x_i}(v')=(a_i + b_i + 1)/2 $ and $ \deg_{x_i}(w')  = (a_i+b_i -1)/2$. Hence, $\deg_{x_i}(v'), \deg_{x_i}(w')  \leq c_i$. This means that $v', w'$ are $\cb$-bounded monomials.

Now, since $vw=v'w' $ and $v,w \in G(B^\cb(u_k))$, it follows from \cite[Lemma~2.7]{De}  that $v',w' \in B^\cb(u_k)$, and since $ v, w, v'$ and $w'$  are in the same degree, we get $v',w' \in G(B^\cb(u_k))$.

Finally,  since $ I=B^{\cb}(u_1,\ldots,u_m)=B^{\cb}(u_1)+\cdots +B^{\cb}(u_m)$, the assertion follows.
 \end{proof}
\begin{Corollary}
The statements of Theorem~\ref{main} hold for equigenerated $\cb$-bounded strongly stable monomial ideals.
\end{Corollary}
\begin{Remark}
\label{reza}
{\em (a) An equigenerated $\cb$-bounded stable ideal is not necessarily sortable. For example, the ideal $I=(x_1^3, x_1^2x_2, x_1x_2^2, x_1x_2x_3) \subset K[x_1,x_2,x_3]$ is  a $\cb$-bounded stable ideal of degree $3$, where $\cb=(3,2,1)$. Note that $\sort (x_1^3,x_1x_2x_3)= (x_1^2x_2,x_1^2x_3) \notin G(I) \times G(I)$.}

{\em (b) Let $v$ be a monomial in $S$  and $I=B(v)$. It followes from \cite[Lemma~2.7]{De} that $G(I)$ is a sortable set. So, since for equigenerated monomials $u_1,\ldots, u_m$ we have $B(u_1, \ldots, u_m)=B(u_1)+\cdots +B(u_m)$,  the statements of the Theorem~\ref{main} hold for equigenerated strongly stable monomial ideals.}
\end{Remark}
 Now we come with an important class of equigenerated $\cb$-bounded strongly stable monomial ideals, called Veronese type ideals.  Let $n$  be a positive integer,  $d$ be an integer, and  $\ab=(a_1,\ldots,a_n)$ be an  integer vector with
 $a_1\geq a_2\geq \cdots\geq a_n$.  The monomial ideal  $I_{\ab,n,d}\subset S=K[x_1,\ldots,x_n]$ with the minimal generating set
\[
G(I_{\ab,n,d})=\{x_1^{b_1}x_2^{b_2}\cdots x_n^{b_n}\; \mid \; \sum_{i=1}^nb_i=d \text{ and  $b_i\leq a_i$ for $i=1,\ldots,n$}\}
\]
 is called a {\em Veronese type} ideal. It is obvious that $I_{\ab,n,d}$ is $\ab$-bounded strongly stable.
\begin{Corollary}
The statements of  Theorem~\ref{main} hold for Veronese type ideals.
\end{Corollary}

\begin{Example}
In the Example~\ref{M}, the ideal $I$ is Veronese type $I_{\ab,3,3}$ with $\ab=(3,2,1)$. Therefore, 
\[
J=( t_{u_1} t_{u_4}- t_{u_2}^2, t_{u_1} t_{u_5}- t_{u_2} t_{u_3},  t_{u_2} t_{u_5}- t_{u_3} t_{u_4},  t_{u_1} t_{u_6}- t_{u_3}t_{u_4},  t_{u_2} t_{u_6}- t_{u_4}t_{u_5},  t_{u_3} t_{u_6}- t_{u_5}^2),
\]
which is confirmed by CoCoA.
\end{Example}

\section{Toric ideals generated by quadratic binomials}
\label{2}
Let $I \subset K[x_1, \ldots, x_n]$ be an equigenerated ideal, such that the toric defining ideal $J$ of $F(I)$ is generated by quadratic binomials. We associate to $I$ a matrix, and show that $J$  is generated by the set of  binomial $2$-minors of this matrix. Indeed,  $J$ is generated by the set of the determinants of $2 \times 2$  submatrices of this matrix which have no zero entries. The construction of the associated matrix  when $n=2$ is different from the cases $n \geq 3$.
For this purpose we introduce the following notation.

\begin{Notation}
\label{mat}
Let $I \subset K[x_1,x_2]$ be an  ideal generated in degree $d$ with the minimal set of monomial generators $G(I)=\{u_1, \ldots, u_{q}\}$  which can be considered as a subset of $G(\mm^d)$. We assume that $I$ contains $x_{1}^d$, because otherwise there exist a positive integer $ d' $  and an ideal $J$ such that  $I=x_{2}^{d'}J$ and  $G(J)$ contains  $x_{1}^{d-d'}$, for which we have $F(I)=F(J)$. Also, we assume that 
\[
u_1=x_{1}^d >_{lex} u_2 = x_1^{d-a}x_2^a >_{lex} u_3 \ldots >_{lex} u_{q-1} >_{lex} u_q=x_1^{d-a-b}x_2^{a+b },
\]
where $1 \leq a,b \leq d-1$ and $2 \leq a+b \leq d$.

We arrange the columns of the matrix $\mathcal{M}$ in the following way
\begin{align*}
\mathcal{M}=\begin{pmatrix}
x_1^d=u_1& x_1^{d-1}x_2 & x_1^{d-2}x_2^2 & \ldots  &  x_1^{d-b}x_2^{b } \\ 
x_1^{d-1}x_2 & x_1^{d-2}x_2^2 & x_1^{d-3}x_2^3  & \ldots  &   x_1^{d-b-1}x_2^{b +1}\\ 
\vdots &\vdots   &\vdots \\
 x_1^{d-a}x_2^a =u_2&   x_1^{d-a-1}x_2^{a+1}  &   x_1^{d-a-2}x_2^{a+2}&\ldots  &  x_1^{d-a-b}x_2^{a+b }=u_q
\end{pmatrix} .
\end{align*}
 We replace the enteries of $\mathcal{M}$  belonging to $G(\mm^d) \setminus G(I)$  by $0$  and denote the obtained matrix by  $\mathcal{M}_I$.  We also replace any nonzero element $u$  of $\mathcal{M}_I$  by the indeterminate $t_u$ of $R=K[t_{u_1}, \ldots, t_{u_q}]$ and denote this matrix by $\mathcal{T}_I$ .
\end{Notation}
\begin{Theorem}
\label{tor}
Let $I \subset K[x_1,x_2]$ be a monomial ideal generated in degree $d$ with the unique minimal set of monomial generators $G(I)=\{u_1, \ldots, u_{q}\}$ and let $F(I)=K[t_{u_1}, \ldots, t_{u_q}]/J$. If the toric ideal $J$ is generated by quadratic binomials, then $J$ is the ideal generated by the set of  binomial $2$-minors of $\mathcal{T}_I$.
\end{Theorem}
\begin{proof}
Let $m_j$ be the least common multiple of all entries of the $j$-th column of $\mathcal{M}$ for  $j=1 , \ldots, b+1$. Dividing the $j$-th column of $\mathcal{M}$ by $m_j$ for all $j$, we get a matrix whose  entries of all columns  are the monomials of $G(\mm^{a})$ ordered lexicographically up to down. So, evey 2-minor of this matrix is zero. This implies that,   every  binomial $2$-minor of $\mathcal{T}_I$  belongs to $J$.

Conversely, we show that any quadratic binomial of $G(J)$ stands as the determinant of a $2 \times 2$ submatrix of $\mathcal{T}_I$ which has no zero entries. Let $f=t_ut_v-t_{u'}t_{v'} \in J $, where $u,v, u',v' \in G(I)$. It is clear from the arrangement of the columns of $\mathcal{M}$ that $u$ and $v$ appear on  the main diagonal of a $2 \times 2$ submatrix of $\mathcal{M}$ (note that $uv \neq x_1^{2d-1}x_2$ and also $u'v' \neq x_1x_2^{2d-1}$). We need to show that $u'$ and $v'$ appear on  the secondary diagonal of  the same submatrix.  Let $u=x_1^{d-p}x_2^p, v=x_1^{d-q}x_2^q, u'=x_1^{d-r}x_2^r, v'=x_1^{d-s}x_2^s$. Since $f \in J$, it follows that $uv=u'v'$ and hence $p+q=r+s$ . So, without loss of generality, we may assume that $s>p$ and $q>r$. In addition, we let $u$ and $v$ be the $ij$-th and the $kl$-th entries of  $\mathcal{M}$ respectively. We show that $u'$ and $v'$ are the $il$-th and the $kj$-th entries of  $\mathcal{M}$ respectively. Since $p+q=r+s$, therefore  $s-p=q-r$. So, it follows from the arrangement of the columns of $\mathcal{M}$  that $u'$ and $v'$ appear as the $il$-th and the $kj$-th entries of  $\mathcal{M}$ respectively, and the proof is complete. 
\end{proof}
In the next theorem we let $I$ be an equigenerated monomial ideal in the polynomial ring $S=K[x_1, \ldots, x_n]$ with $n\geq 3$,  and $T_I$ be its associated matrix introduced in Notation~\ref{associate}.
\begin{Theorem}
\label{more}
Let $I \subset S=K[x_1, \ldots, x_n]$ with $n\geq 3$ be a monomial ideal generated in degree $d$ with the minimal set of monomial generators $G(I)=\{u_1, \ldots, u_{q}\}$, and let $F(I)=K[t_{u_1}, \ldots, t_{u_q}]/J$. If the toric ideal $J$ is generated by quadratic binomials, then $J$ is the ideal generated by the set of  binomial $2$-minors of $T_I$.
\end{Theorem}
\begin{proof}
As we stated in the proof of Theorem~\ref{main}, if we divide the entries of the $i$-th row of the matrix $M$ by $x_i$  for $i=1, \ldots , n$,  we get a matrix whose entries of all rows  are the monomials of $G(\mm^{d-1})$ ordered lexicographically from left to right. This implies that all 2-minors of $M$ are zero and therefore, any binomial $2$-minors of $T_I$ is contained in $J$.

Conversely, for the nonzero monomials $u, v, u', v' \in G(I)$, let the nonzero binomial $f=t_ut_v-t_{u'}t_{v'} $ belongs to $J$. We show that $f$ is a 2-minor of $T_I$. Set $u=x_1^{{\alpha}_1} \ldots x_n^{{\alpha}_n}, v=x_1^{{\beta}_1} \ldots x_n^{{\beta}_n}, u'=x_1^{{\alpha'}_1} \ldots x_n^{{\alpha'}_n}$ and $v'=x_1^{{\beta'}_1} \ldots x_n^{{\beta'}_n}$.  It is clear that $t_ut_v-t_{u'}t_{v'}  \in J$ if and only if $uv=u'v'$. Therefore, $\alpha_i+ \beta_i= \alpha'_i+ \beta'_i $  for $i=1, \ldots , n$.  Since $f$ is a nonzero binomial, then $uv$ is not pure power of a variable and so there are indices $k\neq l$ such that $x_k | u$ and $x_l |v$. Thus, $x_k | u'v'$ and $x_l | u'v'$. We distinguish the following cases:

i) If $x_k | u'$ and $x_l |v'$, then $u, u'$ appear in the $k$-th row of $M$, and $v, v'$ appear in $l$-th row of $M$. We need to show that $v'$ appears  in the same column of  $u$, and  $v$ appears in the same column of $u'$.  Since $\alpha_i+ \beta_i= \alpha'_i+ \beta'_i $ , we get  $\alpha_i-\beta'_i  = \alpha'_i-\beta_i $ for $i=1, \ldots , n$. Set 
$$u/lcm(u,u')=\bar{u}=x_1^{\bar{\alpha}_1} \ldots x_n^{\bar{\alpha}_n},$$ 
 $$u'/lcm(u,u')=\bar{u'}=x_1^{\bar{\alpha'_1}} \ldots x_n^{\bar{\alpha'}_n},$$ 
$$v/lcm(v',v)=\bar{v}=x_1^{\bar{\beta}_1} \ldots x_n^{\bar{\beta}_n},$$
 and 
$$v'/lcm(v',v)=\bar{v'}=x_1^{\bar{\beta'_1}} \ldots x_n^{\bar{\beta'}_n}.$$ 
For $i=1, \ldots , n$, it is clear that  $\bar{\alpha}_i \neq 0$  if and only if $\bar{\alpha'}_i=  0$,  and also  $\bar{\beta'}_i \neq 0$ if and only if  $\bar{\beta}_i =0$. Let $\bar{\alpha}_i \neq 0$.  Then $\bar{\alpha'}_i=  0$. Now,  $\bar{\beta'}_i \neq 0$,  becauce otherwise the equality $\alpha_i-\beta'_i  = \alpha'_i-\beta_i $ gives a contradiction, since the left side is positive and the right side is negative. Therefore, $\beta_i =0$ and so $\alpha_i-\beta'_i  = \alpha'_i-\beta_i =0$. It follows that $\bar{u}=\bar{v'}$ and $\bar{u'}=\bar{v}$. Now, since $u=lcm(u,u')\bar{u}$,  $u'=lcm(u,u')\bar{u'}$ and also $v'=lcm(v',v)\bar{v'}$,  $v=lcm(v',v)\bar{v'}$, it follows that$\begin{pmatrix}
u& u'\\ 
v'& v\\ 
\end{pmatrix}$  is a submatrix of $M$ and hence $f$ is a 2-minor of $T_I$.

ii) Let  $x_k$, $x_l$  do not divide $v'$. It follows that $x_k, x_l |u'$. So, there exists an index $t \neq k,l$ such that $x_t |v'$. Therefore, $x_t| uv$. If $x_t|v$,  then $u, u'$ appear in the $k$-th row of $M$ and $v, v'$  appear in $t$-th row of $M$, and the conclution is exactly the same as in the case (i). Now, assume that $x_t$ does not divide $v$. So $x_t|u$. Therefore,  $u, v'$ appear in the $t$-th row of $M$, and $u', v$ appear in the $l$-th row of $M$. We need to show that $u'$ appears  in the same column of  $u$, and  $v$ appears in the same column of $v'$.  Since $\alpha_i+ \beta_i= \alpha'_i+ \beta'_i $ , we get  $\alpha_i-\alpha'_i  = \beta'_i-\beta_i $ for $i=1, \ldots , n$. Set 
$$u/lcm(u,v')=\hat{u}=x_1^{\hat{\alpha}_1} \ldots x_n^{\hat{\alpha_n}},$$ 
$$u'/lcm(u',v)=\hat{u'}=x_1^{\hat{\alpha'_1}} \ldots x_n^{\hat{\alpha'}_n},$$ 
$$v/lcm(u',v)=\hat{v}=x_1^{\hat{\beta}_1} \ldots x_n^{\hat{\beta}_n},$$ 
and 
$$v'/lcm(v,v')=\hat{v'}=x_1^{\hat{\beta'_1}} \ldots x_n^{\hat{\beta'}_n}.$$
For $i=1, \ldots , n$, it is clear that  $\hat{\alpha}_i \neq 0$  if and only if $\hat{\beta'}_i=  0$,  and also  $\hat{\alpha'}_1 \neq 0$ if and only if  $\hat{\beta}_i =0$. Let $\bar{\alpha}_i \neq 0$.  Then $\hat{\beta'}_i=  0$. Now,  $\hat{\alpha'}_i \neq 0$,  becauce otherwise the equality $\alpha_i-\alpha'_i  = \beta'_i-\beta_i $ gives a contradiction, since the left side is positive and the right side is negative. Thus,  $\beta_i =0$ and hence $\alpha_i-\alpha'_i  = \beta'_i-\beta_i =0$. Therefore, $\hat{u}=\hat{u'}$ and $\hat{v'}=\hat{v}$. Now, since $u=lcm(u,v')\hat{u}$,  $v'=lcm(u,v')\hat{v'}$ and also $u'=lcm(u',v)\hat{u'}$,  $v=lcm(u',v')\hat{v}$, it follows that
$\begin{pmatrix}
u& v'\\ 
u'& v\\ 
\end{pmatrix}$  is a submatrix of $M$ and hence $f$ is a 2-minor of $T_I$.
\end{proof}
An important  consequence of Theorem~\ref{tor} and Theorem~\ref{more} is a characterization of the fiber cone of Freiman ideals.

We recall that the {\em analytic spread} $\ell (I)$ of an ideal $I$  is by definition the Krull dimension of $F(I)$. The following definition is obtained from \cite{HZ1}.
\begin{Definition}
An equigenerated monomial ideal $I$ is called  a  Freiman ideal, if $\mu (I^2)=\ell (I) \mu (I)- {\ell (I) \choose 2}$.
\end{Definition}

\begin{Corollary} 
{\em Assume that  $T_I$  and $\mathcal{T_I}$ are the matrices  introduced in Notation~\ref{associate} and Notation~\ref{mat}.}
\begin{enumerate}
{\em \item[(a)] Let $I=(u_1, \ldots, u_q) \subset K[x_1,x_2]$ be a Freiman ideal with the fiber cone $F(I)=K[t_{u_1}, \ldots, t_{u_q}]/J$. Then, the toric ideal $J$ is generated by the set of  binomial $2$-minors of $\mathcal{T_I}$.}

{\em \item[(b)]   Let $I=(u_1, \ldots, u_q) \subset K[x_1, \ldots, x_n]$ with $n \geq3$ be a Freiman ideal and $F(I)=K[t_{u_1}, \ldots, t_{u_q}]/J$ be its fiber cone. Then, the toric ideal $J$ is generated by the set of  binomial $2$-minors of $T_I$.
}
\end{enumerate}
\end{Corollary}
\begin{proof} (a) , (b). Let $I$ be a Freiman ideal. Then, the toric defining ideal $J$ of $F(I)$ is generated by binomials, (e.g., see \cite[Lemma~5.2]{EH}). On the other hand, $J$ has a 2-linear resolution by \cite[Theorem~2.3]{HHZ}. Therefore,  $J$ is generated by quadratic binomials. Now, (a) follows from Theorem~\ref{tor} and (b) follows from Theorem~\ref{more}. 
\end{proof}
\begin{Remark}
 In \cite{OH}  there exists an example  of a non-Koszul square-free semigroup ring whose toric ideal  is generated by quadratic binomials but possesses no quadratic Gr\"obner basis (\cite[Example 2.1]{OH}). Therefore, the set of  binomial $2$-minors of $T_I$ in Theorem~\ref{more} may not be a Gr\"obner basis of the toric ideal $J$.
\end{Remark}
\begin{Remark}
Theorem~\ref{main} may fails when $I \subset K[x_1,x_2]$ is not generated by  a sortable set of monomials, even the defining ideal of $F(I)$ is generated by quadratic binomials. For example, let $I=(x_{1}^{5}, x_{1}^3x_{2}^2,  x_{1}^2x_{2}^3,  x_{1}x_{2}^4)$. We set  $u_1=x_{1}^{5}, u_2= x_{1}^3x_{2}^2, u_3= x_{1}^2x_{2}^3$ and $ u_4 = x_{1}x_{2}^4$. Note that $G(I)$ is not sortable, since $\sort(u_1,u_2)= (x_{1}^4x_{2}, x_{1}^4x_{2}) \notin G(I)\times G(I)$. One can check by CoCoA that $F(I)= K[t_{u_1}, \ldots, t_{u_4}]/ J$ where $J=(t_{u_1}t_{u_4}- t_{u_2}^2, t_{u_2}t_{u_4}- t_{u_3}^2)$. While,
\[
T_I=\begin{pmatrix}
t_{u_1}& 0 & t_{u_2} & t_{u_3} &t_{u_4}\\ 
0 & t_{u_2}   & t_{u_3} & t_{u_4}& 0\\ 
\end{pmatrix}.
\]
So, the ideal generated by the set of  binomial $2$-minors of  $ T_I$  is $(t_{u_2}t_{u_4}- t_{u_3}^2)$. Thus, in this case we use Theorem~\ref{tor} to find $J$. We have 
\[
\mathcal{T}_I=\begin{pmatrix}
t_{u_1} & 0 & t_{u_2} \\
0 & t_{u_2}   & t_{u_3} \\
t_{u_2} & t_{u_3}   & t_{u_4} \\
\end{pmatrix}.
\]
It is easy to see that $J$ is the ideal generated by the set of  binomial $2$-minors of $ \mathcal{T}_I$.
\end{Remark}
\begin{Example}
{\em (a) Let  $I=(x_{1}^{12}, x_{1}^9x_{2}^3,  x_{1}^6x_{2}^6,  x_{1}^3x_{2}^9)\subset K[x_1,x_2]$ . The ideal $I$ is a Freiman ideal, since $\ell(I)=2$ and $\mu(I^2)=7=2\mu (I)- {2\choose 2}$.   Set  $u_1=x_{1}^{12}, u_2= x_{1}^9x_{2}^3, u_3= x_{1}^6x_{2}^6$ and $ u_4= x_{1}^3x_{2}^9$ . Checking by CoCoA we get $F(I)= K[t_{u_1},\ldots, t_{u_4}]/ J$ where 
\[
J=(t_{u_1}t_{u_3}- t_{u_2}^2, t_{u_1}t_{u_4}- t_{u_2}t_{u_3}, t_{u_2}t_{u_4}- t_{u_3}^2). 
\]
Note that $G(I)$ is not sortable, since $\sort(u_1,u_2)= (x_{1}^{11}x_{2}, x_{1}^{10}x_{2}^2) \notin G(I)\times G(I)$. Using Theorem~\ref{tor} to find $J$ we get
\[
\mathcal{T}_I=\begin{pmatrix}
t_{u_1}& 0 & 0 &  t_{u_2}& 0 & 0 &  t_{u_3}\\ 
0 & 0  &t_{ u_2  }&  0 & 0 & t_{u_3 }& 0\\ 
0 & t_{u_2}   & 0 & 0 & t_{u_3} & 0 &  0\\ 
t_{u_2} & 0 & 0   & t_{u_3} & 0 & 0 &  t_{u_4}\\ 
\end{pmatrix}.
\]
We see that $J$ is the ideal generated by the set of  binomial $2$-minors of $ \mathcal{T}_I$.
}

{\em (b)
  Let  $I=(x_{1}^{3}, x_{1}^2x_3,  x_{1}x_{3}^2, x_{2}^3,  x_{3}^3)\subset K[x_1, x_2 ,x_3]$ . Then $I$ is a Freiman ideal, because $\ell(I)=3$ and $\mu(I^2)=12=3\mu (I)- {3 \choose 2}$.  Set  $u_1=x_{1}^{3}, $ $u_2= x_{1}^2x_{3}, $ $ u_3= x_{1}x_{3}^2$, $  u_4= x_{2}^3$ and $ u_5 = x_{3}^3$ . Let $F(I)=K[u_1, \ldots u_5]/J$.  It follows from  Theorem~\ref{more} that  $J$ is the ideal generated by the set of  binomial $2$-minors of $ T_I$,  where
\[
T_I=\begin{pmatrix}
t_{u_1}& 0 & t_{u_2}&  0 & 0 &  t_{ u_3}\\ 
0 & 0  & 0 &t_{ u_4}& 0 & 0\\ 
 t_{ u_2}& 0   & t_{ u_3} & 0 & 0&  t_{u_5}\\ 
\end{pmatrix}.
\]
So $J=(t_{u_1}t_{u_3}- t_{u_2}^2, t_{u_1}t_{u_5}- t_{u_2}t_{u_3}, t_{u_2}t_{u_5}- t_{u_3}^2).$ The result is confirmed by CoCoA. Note that $G(I)$ is not sortable, since $\sort(u_1,u_4)= (x_{1}^{2}x_{2}, x_{1}x_{2}^2) \notin G(I)\times G(I)$.
}
\end{Example}

\end{document}